\documentclass[11pt]{article}
\usepackage{amsmath,amssymb,amscd,enumerate,eucal,bm}
\usepackage{color}
\numberwithin{equation}{section}
\newcommand{\BH}{\mathbb{H}}

\begin{document}
\setlength{\baselineskip}{15pt}
\title{A Note on Commutators of the Fractional Sub-Laplacian on Carnot Groups}
\author{Ali Maalaoui$^{1}$ }
\addtocounter{footnote}{1}
\footnotetext{Department of mathematics and natural sciences, American University of Ras Al Khaimah, PO Box 10021, Ras Al Khaimah, UAE. E-mail address:
{\tt{ali.maalaoui@aurak.ac.ae}}}

\date{}
\maketitle
\begin{abstract}
In this manuscript, we provide a point-wise estimate for the $3$-commutators involving fractional powers of the sub-Laplacian on Carnot groups of homogeneous dimension $Q$. This can be seen as a fractional Leibniz rule in the sub-elliptic setting. As a corollary of the point-wise estimate, we provide an $(L^{p},L^{q})\to L^{r}$ estimate for the commutator, provided that $\frac{1}{r}=\frac{1}{p}+\frac{1}{q}-\frac{\alpha}{Q}$.
{\footnotesize \baselineskip 4mm }
\end{abstract}

\medskip

{\footnotesize\textit{}}

\newtheorem{theorem}{Theorem}[section]
\newtheorem{lemma}{Lemma}[section]
\newtheorem{corollary}{Corollary}[section]
\newtheorem{definition}{Definition}[section]
\newtheorem{proposition}{Proposition}[section]
\newtheorem{remark}{Remark}[section]
\newtheorem{claim}{Claim}[section]

\section{Introduction and Main Results}
In this paper we propose to investigate commutator type estimates for fractional powers of the sub-Laplacian on a Carnot group $\mathbb{G}$. In fact, we will be focusing mainly on a fractional type Leibniz rule. These kind of estimates are very effective in studying the regularity or certain fractional PDEs. The study of these type of PDEs can be of many interests, namely, from an analytical point of view, but also from a geometric point of view, after the paper of \cite{GG}, where a family of sub-elliptic fractional conformally invariant operators were exhibited. These operators are the CR parallel to the GJMS operators \cite{GJMS} defined in the Riemannian setting.

Because of this geometric motivation, the Riemannian GJMS operators where intensively investigated and in the special case of the euclidean space $\mathbb{R}^{n}$, the fractional powers of the Laplacian became a main research focus. For instance, fractional type Leibniz rules were proved and investigated in \cite{Schik} using potential analysis and then another alternative proof was provided using harmonic extensions in \cite{lenz}. 
As is the case of local operators, the main tool to study regularity of certain PDEs is localizing the functions with a cut-off. For non-local operators this procedure becomes more complicated. That is why a Leibnitz type rule is valuable in this case. In our work, we will follow the potential analysis approach developed in \cite{Schik} in order to provide a fractional sub-elliptic Leibniz type rule in a general Carnot group. Then, in Section 6 we discuss the case of the Heisenberg group and the two kinds of fractional type operators related to the sub-Laplacian.\\

We start here by defining the commutators and introducing the main results. Given a Carnot group $\mathbb{G}$ we denote by $Q$ its homogeneous dimension. For this purpose, we fix $0<\alpha<Q$, $u,v\in S(\mathbb{G})$, the 3-commutator $H_{\alpha}(\cdot,\cdot)$ is then defined by
$$H_{\alpha}(u,v)=(\Delta_{b})^{\frac{\alpha}{2}}(uv)-u(\Delta_{b})^{\frac{\alpha}{2}}v-v(\Delta_{b})^{\frac{\alpha}{2}}u.$$
 Also, given $0<\tau<Q$ and $\beta,\delta>0$ such that $\beta+\alpha<\min\{1,\tau\}$ we define the operator $T_{\tau,\beta,\delta}$ by

\begin{align}
T_{\tau,\beta,\delta}(u,v)&=[(\Delta_{b})^{-\frac{\tau}{2}}u,(\Delta_{b})^{\frac{\beta}{2}}](\Delta_{b})^{\frac{\delta}{2}}v\notag\\
&=(\Delta_{b})^{-\frac{\tau}{2}}u(\Delta_{b})^{\frac{\beta+\delta}{2}}v-(\Delta_{b})^{\frac{\beta}{2}}((\Delta_{b})^{-\frac{\tau}{2}}u(\Delta_{b})^{\frac{\delta}{2}}v)
\end{align}

Our first result, involves a point-wise estimate, similar to the one in \cite{Schik}, that can be seen as a fractional type Leibniz rule for $0<\alpha<2$.
\begin{theorem}
Assume that $u$ and $v$ are two Schwartz functions on $\mathbb{G}$ and let $0<\alpha<2$ and $\epsilon>0$. Given $\tau_{1}$ and $\tau_{2}$ in $(\max\{0,\alpha-1\},\alpha]$ such that $\tau_{1}+\tau_{2}>\alpha$, there exists $L\in \mathbb{N}$, $s_{j,1}\in (0,\tau_{1})$, $s_{j,2}\in (0,\tau_{2})$, for $j=1,\cdots, L$, satisfying $\tau_{1}+\tau_{2}-s_{j,1}-s_{j,2}-\alpha\in [0,\epsilon)$ such that
\begin{equation}
|H(u,v)|(x)\lesssim \sum_{j=1}^{L}R_{\tau_{1}+\tau_{2}-s_{j,1}-s_{j,2}-\alpha}\Big(R_{s_{j,1}}|(\Delta_{b})^{\frac{\tau_{1}}{2}}u|R_{s_{j,2}}|(\Delta_{b})^{\frac{\tau_{2}}{2}}v|\Big)(x).
\end{equation}

\end{theorem}

One of the tools that allow us to extend Theorem 1.1 to higher values of $\alpha$, is the following commutator estimates involving the Riesz potential. This type of estimates can be compared to the result of Chanillo \cite{Chan}, though, the estimates that were provided there for the Euclidean case are sharp.
\begin{theorem}
Assume that $u$ and $v$ are two Schwartz functions. Given $\tau>0$ and $\beta$ and $\delta$ two non-negative numbers such that $\beta+\delta<\min\{\tau,1\}$. Given $\epsilon>0$, there exists $L\in \mathbb{N}$, $s_{j,1},\tilde{s}_{j,1}>0$ and $s_{j,2},\tilde{s}_{j,2}\in (0,\tau)$ satisfying $\tilde{s}_{j,1}+\tilde{s}_{j,2}=s_{j,1}+s_{j,2}=\tau-\beta-\delta$ and $\tilde{s}_{j,1}<\epsilon$, such that
\begin{equation}
|T_{\tau,\beta,\delta}(u,v)|(x)\lesssim\sum_{j=1}^{L}R_{s_{j,1}}|u|(x)R_{s_{j,2}}|v|(x)+R_{\tilde{s}_{j,1}}(|v|R_{\tilde{s}_{j,1}}|u|)(x).
\end{equation}
\end{theorem}
A combination of Theorem 1.1 ans 1.2 provides the following
\begin{corollary}
Under the assumptions of Theorem 1.1, given $0<\alpha<Q$ and $\epsilon>0$. Given $\tau_{1}$ and $\tau_{2}$ in $(\max\{0,\alpha-1\},\alpha]$ such that $\tau_{1}+\tau_{2}>\alpha$, there exists $L\in \mathbb{N}$, $s_{j,1}\in (0,\tau_{1})$, $s_{j,2}\in (0,\tau_{2})$, for $j=1,\cdots, L$, satisfying $\tau_{1}+\tau_{2}-s_{j,1}-s_{j,2}-\alpha\in [0,\epsilon)$ such that
\begin{equation}
|H(u,v)|(x)\lesssim \sum_{j=1}^{L}R_{\tau_{1}+\tau_{2}-s_{j,1}-s_{j,2}-\alpha}\Big(R_{s_{j,1}}|(\Delta_{b})^{\frac{\tau_{1}}{2}}u|R_{s_{j,2}}|(\Delta_{b})^{\frac{\tau_{2}}{2}}v|\Big)(x).
\end{equation}
\end{corollary}

Using the $\lambda$-kernel estimates in \cite{FS}, we easily show that
\begin{corollary}
Let $u$ and $v$ be two functions in the Schwartz class of $\mathbb{G}$, and $\alpha \in (0,Q)$, then
\begin{equation}\label{intest}
\|H_{\alpha}(u,v)\|_{p}\lesssim \|(\Delta_{b})^{\frac{\alpha}{2}}u\|_{q_{1}}\|(\Delta_{b})^{\frac{\alpha}{2}}v\|_{q_{2}}
\end{equation}
for $\frac{1}{p}=\frac{1}{q_{1}}+\frac{1}{q_{2}}-\frac{\alpha}{Q}$.
\end{corollary}

\section{Preliminaries and Notations}
In what follows, $\mathbb{G}$ will be a Carnot group with homogeneous dimension $Q$ and $\Delta_{b}$ will be its sub-Laplacian as defined in Folland \cite{Fol1} and Folland-Stein in \cite{FS}. We will start by recalling some definitions and properties of the fractional sub-Laplacian and adapt the following comparaison notations:
We write $f(x)\approx g(x)$, if there exists a constant $C>0$ such that
$$\frac{1}{C}f(x)\leq g(x)\leq Cf(x).$$
Also, we will write $f(x)\lesssim g(x)$ if there exists a constant $C>0$ such that $$f(x)\leq Cg(x).$$
Recall that there exists a gauge $|\cdot|$ on $\mathbb{G}$ that we can assume symmetric, that induces a quasi-distance (see \cite{Bon} Chap 5, for more details). Moreover we have that the existence of two positive constants $c$ and $C$ such that
$$c||x|-|y||\leq |yx|\leq C(|x|+|y|),$$
We can assume without loss of generality that $c<1$ and $C>1$. We provide now a very usefull lemma, that we will be using later on in our computations
\begin{lemma}[\cite{Fol1}]
Let $f$ be a $\lambda$-homogeneous function, that is smooth on $\mathbb{G}\setminus\{0\}$ then there exists a constant $C>0$ such that if $|xy|\approx |x|$ we have
$$|f(xy)-f(x)|\leq C\max\{|xy|^{\lambda-1},|x|^{\lambda-1}\}|y|$$
\end{lemma}
{\it Proof:}
A version of this Lemma was proved in \cite{Fol1} for $2|y|<|x|$. In fact, using the homogeneity of $f$, we can rescale and assume that $|x|=1$ and $|y|\leq \frac{1}{2}$, therefore $|xy|\not=0$ and hence, the proof follows from the regular intermediate calue theorem, leading to
$$|f(xy)-f(x)|\leq C|x|^{\lambda-1}|y|.$$
 It remains then to prove the case when $|y|\approx |x|$. Again, if $|xy|>\frac{1}{4}$, then the proof is done, or else, we have that $|xy|<\frac{1}{4}$, then $|y|>\frac{1}{2}$, hence, we have that
$$|f(xy)-f(x)|\leq C|xy|^{\lambda-1}|y|.$$
\hfill $\Box$

We move now to the definition of the fractional sub-Laplacian and its inverse. For this purpose, we consider the heat semi-group defined by $H_{t}=e^{t\Delta_{b}}$, then one can define as in \cite{Fol1,FS}, the following fractional operators for $\alpha>0$:
$$(\Delta_{b})^{-\alpha}=\frac{1}{\Gamma(\alpha)}\int_{0}^{\infty}t^{\alpha-1}H_{t}dt$$
and
$$(\Delta_{b})^{\alpha}=\frac{1}{\Gamma(\alpha)}\int_{0}^{\infty}t^{\alpha-k-1}(\Delta_{b})^{k}H_{t}dt,$$
for any integer $k>\alpha$.
We have now the following:
\begin{proposition}[\cite{Fol1,FS}]
Let $0<\alpha<Q$ and consider $h(t,x)$ the fundamental solution of the operator $\Delta_{b}+\frac{\partial}{\partial t}$, then the integral
$$R_{\alpha}(x)=\frac{1}{\Gamma(\frac{\alpha}{2})}\int_{0}^{\infty}t^{\frac{\alpha}{2}-1}h(t,x)dt$$
converges absolutely and it satisfies the following properties:
\begin{itemize}
\item $R_{\alpha}$ is a kernel of type $\alpha$. In particular it is homogeneous of degree $\alpha-Q$
\item $R_{2}$ is the fundamental solution of $\Delta_{b}$
\item $R_{\alpha}*R_{\beta}=R_{\alpha+\beta}$ for $\alpha$ and $\beta>0$.
\item For $f\in L^{p}(\mathbb{G})$ and $1<p<\infty$, we have that $$(\Delta_{b})^{-\frac{\alpha}{2}}f=f*R_{\alpha}.$$
\end{itemize}
\end{proposition}
As a corollary of this proposition one has
$$R_{\alpha}(x)\approx |x|^{-Q+\alpha},$$
and $\rho(x)=(R_{\alpha}(x))^{\frac{1}{\alpha-Q}}$ defines a $\mathbb{G}$-homogeneous norm, smooth away from the origin and it induces a quasi-distance that is equivalent to the left-invariant Carnot-Caratheodory distance.
In a similar way, one can define the function $\tilde{R}_{\alpha}$, introduced in \cite{Fran}, for $\alpha<0$ and$\alpha \not \in \{0,-2,-4,\cdots\}$ by
$$\tilde{R}_{\alpha}(x)=\frac{\frac{\alpha}{2}}{\Gamma(\frac{\alpha}{2})}\int_{0}^{\infty}t^{\frac{\beta}{2}-1}h(t,x)dt.$$
Again, it is easy to see that $\tilde{R}_{\alpha}$ is $\mathbb{G}$-homogeneous of degree $\alpha-Q$ and
$$\tilde{R}_{\alpha}(x)\approx |x|^{\alpha-Q}.$$
Using this function, one can define another representation for the fractional sub-Laplacian, that will be fundamental in the proof of our results.
\begin{theorem}[\cite{Fran}]
If $u$ is a Schwartz function on $\mathbb{G}$, then for $0<\alpha<2$ one has
$$(\Delta_{b})^{\frac{\alpha}{2}}u(x)=PV\int_{\mathbb{G}}(u(y)-u(x))\tilde{R}_{-\alpha}(y^{-1}x)dy$$
\end{theorem}
For the sake of notation, we will use in this paper the following notation:
$$R_{\alpha}u:=(\Delta_{b})^{-\frac{\alpha}{2}}u=u*R_{\alpha}.$$

\section{Proof of Theorem 1.1 }
We are interested in the study of the commutator
$$H_{\alpha}(u,v)=(\Delta_{b})^{\frac{\alpha}{2}}(uv)-u(\Delta_{b})^{\frac{\alpha}{2}}v-v(\Delta_{b})^{\frac{\alpha}{2}}u.$$
Notice that 
$$H_{\alpha}(u,v)(x)=\int_{\mathbb{G}}\Big[u(x)-u(y)][v(x)-v(y)\Big]\tilde{R}_{-\alpha}(y^{-1}x)dy.$$
We take now $u=(\Delta_{b})^{-\frac{\tau_{1}}{2}}a=R_{\tau_{1}}a$ and $v=(\Delta_{b})^{-\frac{\tau_{2}}{2}}b=R_{\tau_{2}}b$, then we can write
$$H_{\alpha}(u,v)(x)=\int\int\int \Big[R_{\tau_{1}}(\eta^{-1}y)-R_{\tau_{1}}(\eta^{-1}x)\Big]\Big[(R_{\tau_{2}}(\xi^{-1}y)-R_{\tau_{2}}(\xi^{-1}x)\Big]\tilde{R}_{\alpha}(y^{-1}x)a(\eta)b(\xi)d\eta d\xi dy.$$
Therefore, one needs to study the kernel $k$ defined by:
$$k(x,y,\eta,\xi):=\Big[R_{\tau_{1}}(\eta^{-1}y)-R_{\tau_{1}}(\eta^{-1}x)\Big]\Big[(R_{\tau_{2}}(\xi^{-1}y)-R_{\tau_{2}}(\xi^{-1}x)\Big]\tilde{R}_{\alpha}(y^{-1}x).$$
In order to simplify the computations, we introduce the kernel $k_{\tau}$ defined by
$$k_{\tau}(x,y,\eta)=|R_{\tau}(\eta^{-1}y)-R_{\tau}(\eta^{-1}x)|,$$
and we split the space into three parts using the following characteristic functions:
$$\left\{\begin{array}{lll}
\chi_{1}=\chi_{|y^{-1}x|<2|\eta^{-1}y|}\chi_{|y^{-1}x|<2|\eta^{-1}x|}\\
\chi_{2}=\chi_{|y^{-1}x|<2|\eta^{-1}y|}\chi_{|y^{-1}x|>2|\eta^{-1}x|}\\
\chi_{3}=\chi_{|y^{-1}x|>2|\eta^{-1}y|}\chi_{|y^{-1}x|\leq 2|\eta^{-1}x|}
\end{array}
\right.
.$$
We first notice that $$|\eta^{-1}y|\chi_{1}\approx |\eta^{-1}x|\chi_{1}.$$
Therefore, using Lemma 2.1, we have, for $0<\delta<\min\{\tau,1\}$, that
\begin{align}
k_{\tau}(x,y,\eta)\chi_{1}&\lesssim\min\{|\eta^{-1}x|,|\eta^{-1}y|\}^{\tau-Q}|y^{-1}x|\chi_{1}\notag\\
&\lesssim|\eta^{-1}x|^{\tau-Q-\delta}|y^{-1}x|^{\delta}\chi_{1}\notag\\
&\lesssim|\eta^{-1}y|^{\tau-Q-\delta}|y^{-1}x|^{\delta}\chi_{1}.
\end{align}
Moreover,
$$k_{\tau}(x,y,\eta)\chi_{2}\lesssim|\eta^{-1}x|^{\tau-Q}\chi_{2}\lesssim|\eta^{-1}x|^{\tau-Q-\delta}|y^{-1}x|^{\delta}\chi_{2},$$
and similarly 
$$k_{\tau}(x,y,\eta)\chi_{3}\lesssim|\eta^{-1}y|^{\tau-Q}\chi_{3}\lesssim|\eta^{-1}y|^{\tau-Q-\delta}|y^{-1}x|^{\delta}\chi_{3}.$$
Therefore, since $k(x,y,\eta,\xi)=k_{\tau_{1}}(x,y,\eta)k_{\tau_{2}}(x,y,\xi)\tilde{R}_{-\alpha}(y^{-1}x)$ we have, for $0<\delta_{1}<\min\{\tau_{1},1\}$ and $0<\delta_{2}<\min\{\tau_{2},1\}$, 
\begin{align}
k(x,y,\eta,\xi)\lesssim&  |y^{-1}x|^{-Q-\alpha+\delta_{1}+\delta_{2}}\Bigg(|\eta^{-1}x|^{-Q+\tau_{1}-\delta_{1}}|\xi^{-1}y|^{-Q+\tau_{2}-\delta_{2}}\\
&+|\eta^{-1}y|^{-Q+\tau_{1}-\delta_{1}}|\xi^{-1}x|^{-Q+\tau_{2}-\delta_{2}}+|\eta^{-1}y|^{-Q+\tau_{1}-\delta_{1}}|\xi^{-1}y|^{-Q+\tau_{2}-\delta_{2}}\\
&+|y^{-1}x|^{-\delta_{1}-\delta_{2}}|\eta^{-1}x|^{-Q+\tau_{1}}|\xi^{-1}x|^{-Q+\tau_{2}}\chi_{|y^{-1}x|>2\max\{|\eta^{-1}x|;|\xi^{-1}x|\}} \Bigg).
\end{align}
Now notice that using "polar coordinates", if one takes $0<\delta_{1}+\delta_{2}-\alpha<\varepsilon$, where $\varepsilon$ is positive and small, we have for $\epsilon>0$ that
$$\int_{|y^{-1}x|>\epsilon} |y^{-1}x|^{-Q-\alpha}\chi_{|y^{-1}x|>2\max\{|\eta^{-1}x|;|\xi^{-1}x|\}}dy\lesssim\max\{|\eta^{-1}x|,|\xi^{-1}x|\}^{-\alpha}\lesssim|\eta^{-1}x|^{-\delta_{1}}|\xi^{-1}x|^{-\alpha+\delta_{1}}.$$
This yields, 
\begin{align}
|H_{\alpha}(u,v)(x)|&\lesssim \int \int \int k(x,y,\eta,\xi)|a(\eta)||b(\xi)|d\eta d\xi dy\notag\\
&\lesssim\Big(R_{\delta_{1}+\delta_{2}-\alpha}(R_{\tau_{2}-\delta_{2}}|b|)(x)R_{\tau_{1}-\delta_{1}}|a|(x)\notag\\
&\qquad+R_{\delta_{1}+\delta_{2}-\alpha}(R_{\tau_{1}-\delta_{1}}|a|)(x)R_{\tau_{2}-\delta_{1}}|b|(x)\notag\\
&\qquad+R_{\delta_{1}+\delta_{2}-\alpha}(R_{\tau_{1}-\delta_{1}}|a|R_{\tau_{2}-\delta_{2}}|b|)(x)\notag\\
&\qquad+R_{\tau_{1}-\delta_{1}}|a|(x)R_{\tau_{2}+\delta_{1}-\alpha}|b|(x)\Big).
\end{align}
Hence
\begin{align}
|H_{\alpha}(u,v)(x)|\lesssim& R_{\tau_{1}-\delta_{1}}|a|(x)R_{\tau_{2}+\delta_{1}-\alpha}|b|(x)+R_{\tau_{1}+\delta_{2}-\alpha}|a|(x)R_{\tau_{2}-\delta_{1}}|b|(x)\notag\\
&\qquad\qquad+R_{\delta_{1}+\delta_{2}-\alpha}(R_{\tau_{1}-\delta_{1}}|a|R_{\tau_{2}-\delta_{2}}|b|)(x).
\end{align}
Therefore, one can state that
$$|H(u,v)(x)|\leq C\sum_{j=1}^{L}R_{\tau_{1}+\tau_{2}-s_{j,1}-s_{j,2}-\alpha}(R_{s_{1,j}}|a|R_{s_{j,2}}|b|)(x),$$
proving the desired result.

\section{Proof of Theorem 1.2}
From now on, we will let $\tau>0$ and two non-negative numbers $\beta$ and $\delta$ so that $\beta+\delta<\min\{\tau,1\}$. We then, define the commutator $T_{\tau,\beta,\delta}(u,v)$ by
\begin{align}
T_{\tau,\beta,\delta}(u,v)&=\Big[R_{\tau}u,(\Delta_{b})^{\frac{\beta}{2}}\Big](\Delta_{b})^{\frac{\delta}{2}}v\\
&=R_{\tau}u(\Delta_{b})^{\frac{\beta+\delta}{2}}v-(\Delta_{b})^{\frac{\beta}{2}}(R_{\tau}u(\Delta_{b})^{\frac{\delta}{2}}v).
\end{align}
Notice that
\begin{align}
(\Delta_{b})^{\frac{\beta}{2}}(R_{\tau}u(\Delta_{b})^{\frac{\delta}{2}}v)&=\int \Big[R_{\tau}u(x)(\Delta_{b})^{\frac{\delta}{2}}v(x)-R_{\tau}u(y)(\Delta_{b})^{\frac{\delta}{2}}v(y)\Big]\tilde{R}_{-\beta}(y^{-1}x)dy\notag\\
&=R_{\tau}u(x)\int \Big[(\Delta_{b})^{\frac{\delta}{2}}v(x)-(\Delta_{b})^{\frac{\delta}{2}}v(y)\Big]\tilde{R}_{-\beta}(y^{-1}x)dy\notag\\
&\qquad+ \int [R_{\tau}u(x)-R_{\tau}u(y)](\Delta_{b})^{\frac{\delta}{2}}v(y)\tilde{R}_{-\beta}(y^{-1}x)dy.
\end{align}
Therefore
\begin{align}
&T_{\tau,\beta,\delta}(u,v)=\int \Big [R_{\tau}u(x)-R_{\tau}u(y)\Big](\Delta_{b})^{\frac{\delta}{2}}v(y)\tilde{R}_{-\beta}(y^{-1}x)dy\notag\\
&=\int \int \Big [R_{\tau}u(x)-R_{\tau}u(y)\Big ]\Big [v(wy)-v(y)\Big]\tilde{R}_{-\beta}(y^{-1}x)\tilde{R}_{-\delta}(w)dydw\notag \\
&=\int \int \Bigg(\Big[R_{\tau}u(x)-R_{\tau}u(w^{-1}y)\Big]\tilde{R}_{-\beta}(wy^{-1}x)-\Big[R_{\tau}u(x)-R_{\tau}u(y)\Big]\tilde{R}_{-\beta}(y^{-1}x)\Bigg)v(y)\tilde{R}_{-\delta}(w)dw dy\notag\\
&=\int \int \int \Bigg(\Big [R_{\tau}(z^{-1}x)-R_{\tau}(z^{-1}w^{-1}y)\Big]\tilde{R}_{-\beta}(wy^{-1}x)\notag\\
&\qquad \qquad\qquad-\Big[R_{\tau}(z^{-1}x)-R_{\tau}(z^{-1}y)\Big]\tilde{R}_{-\beta}(y^{-1}x)\Bigg)u(z)v(y)\tilde{R}_{-\delta}(w)dwdzdy.
\end{align}
As in the previous proof, we propose to study the kernel $k$ defined by
$$k(x,y,z,w)=\tilde{R}_{-\delta}(w)\Bigg(\Big [R_{\tau}(z^{-1}x)-R_{\tau}(z^{-1}w^{-1}y)\Big]\tilde{R}_{-\beta}(wy^{-1}x)-\Big[R_{\tau}(z^{-1}x)-R_{\tau}(z^{-1}y)\Big]\tilde{R}_{-\beta}(y^{-1}x)\Bigg).$$
But, in this situation, we need a better splitting of the space, adapted to the different quantities involved. So we define the sets
$$A_{1}=\{(x,y,z,w)\in (\BH^{n})^{4};|y^{-1}x|\leq 2|y^{-1}z|\text{ and } |y^{-1}x|\leq 2|z^{-1}x|\},$$
$$A_{2}=\{(x,y,z,w)\in (\BH^{n})^{4};|y^{-1}x|\leq 2|y^{-1}z|\text{ and } |y^{-1}x|> 2|z^{-1}x|\},$$
$$A_{3}=\{(x,y,z,w)\in (\BH^{n})^{4};|y^{-1}x|> 2|y^{-1}z|\text{ and } |y^{-1}x|\leq 2|z^{-1}x|\},$$
and
$$B_{1}=\{(x,y,z,w)\in (\BH^{n})^{4};|wy^{-1}x|\leq 2|wy^{-1}z|\text{ and } |wy^{-1}x|>2|z^{-1}x|\},$$
$$B_{2}=\{(x,y,z,w)\in (\BH^{n})^{4};|wy^{-1}x|\leq 2|wy^{-1}z|\text{ and } |wy^{-1}x|\leq 2|z^{-1}x|\},$$
$$B_{3}=\{(x,y,z,w)\in (\BH^{n})^{4};|wy^{-1}x|> 2|wy^{-1}z|\text{ and } |wy^{-1}x|\leq 2|z^{-1}x|\},$$
and
$$C_{1}=\{(x,y,z,w)\in (\BH^{n})^{4};|w|\geq 4|y^{-1}x|\},$$
$$C_{2}=\{(x,y,z,w)\in (\BH^{n})^{4};\frac{|y^{-1}x|}{4}\leq |w| \leq 4|y^{-1}x|\},$$
$$C_{3}=\{(x,y,z,w)\in (\BH^{n})^{4};4|w| < |y^{-1}x|\}.$$
We can verify that each collection of sets exhaust $(\mathbb{G})^{4}$. In a first step, we will split the kernel $k$ so that
$k\leq k_{1}+k_{2}$ where
$$k_{1}(x,y,z,w)=\tilde{R}_{-\delta}(w)\Big|R_{\tau}(z^{-1}x)-R_{\tau}(z^{-1}w^{-1}y)\Big|\tilde{R}_{-\beta}(wy^{-1}x),$$
and
$$k_{2}(x,y,z,w)=\tilde{R}_{-\delta}(w)\Big|R_{\tau}(z^{-1}x)-R_{\tau}(z^{-1}y)\Big|\tilde{R}_{-\beta}(y^{-1}x).$$
\\
\textbf{Step I}: Estimate on $C_{1}\cup C_{2}$.\\
Notice that on $C_{1}\cup C_{2}$ we have that
$$\int \tilde{R}_{-\delta}(w)\chi_{C_{1}\cup C_{2}}dw\lesssim \int |w|^{-Q+\delta}\chi_{C_{1}\cup C_{2}}dw \lesssim |y^{-1}x|^{-\delta}.$$
Therefore,
$$\int k_{2}\chi_{C^{1}\cup C_{2}}dw \lesssim|R_{\tau}(z^{-1}x)-R_{\tau}(z^{-1}y)||y^{-1}x|^{-Q-\beta-\delta}.$$
Using the same procedure that was done in Section 3, we have the existence of $L>0$ and numbers $s_{k},\tilde{s}_{k}>0$ and $ t_{k},\tilde{t}_{k} \in (0,\tau)$ with $\tilde{s}_{k}<\varepsilon$, and $s_{k}+t_{k}=\tilde{s}_{k}+\tilde{t}_{k}=\tau-\beta-\delta$ so that
$$\int \int \int \chi_{C_{1}\cup C_{2}}k_{2}(x,y,z,w)u(z)v(w)dwdydz\lesssim \sum_{i=1}^{L}R_{s_{i}}|u|(x)R_{t_{i}}|v|(x)+R_{\tilde{s}_{i}}(R_{\tilde{t}_{i}}|u||v|)(x),
$$
which is the required estimate. On the other hand, the estimate for $k_{1}$ is not as straight forward as it is for $k_{2}$. Indeed, we will split the set $C_{1}\cup C_{2}$ using the sets $B_{j}$, $j=1,2,3$.\\
For instance, in the set $(C_{1}\cup C_{2})\cap B_{1}$, we have that
$$|wy^{-1}z|\chi_{(C_{1}\cup C_{2})\cap B_{1}}\leq C(|z^{-1}x|+|wy^{-1}x|)\chi_{(C_{1}\cup C_{2})\cap B_{1}}\leq C|z^{-1}x|\chi_{(C_{1}\cup C_{2})\cap B_{1}},$$
and $$2|wy^{-1}z|\chi_{(C_{1}\cup C_{2})\cap B_{1}}\geq |wy^{-1}x|\chi_{(C_{1}\cup C_{2})\cap B_{1}}\geq c(|z^{-1}x|-|yw^{-1}z|)\chi_{(C_{1}\cup C_{2})\cap B_{1}}.$$
Therefore, $|z^{-1}x|$ and $|wy^{-1}z|$ are comparable and we can use Lemma 2.1, to get for any $\varepsilon \in [0,1]$ 
$$k_{1}(x,y,z,w)\chi_{(C_{1}\cup C_{2})\cap B_{1}}\lesssim |w|^{-Q-\delta}\max\{|z^{-1}x|;|wy^{-1}z|\}^{-Q+\tau-\varepsilon}|wy^{-1}x|^{-Q-\beta+\varepsilon}\chi_{C_{1}\cup C_{2}\cap B_{1}}.$$
Now notice that on $C_{1}$, we have that
$$|y^{-1}x|\chi_{C{1}}\lesssim |wy^{-1}x|\chi_{C_{1}},$$
and on $C_{2}$,
$$|w|\chi_{C_{2}}\approx |y^{-1}x|\chi_{C_{2}}.$$
Thus,
\begin{align}
k_{1}(x,y,z,w)\chi_{(C_{1}\cup C_{2})\cap B_{1}}&\lesssim|w|^{-Q-\delta}|z^{-1}x|^{-Q+\tau-\varepsilon}|y^{-1}x|^{-Q-\beta+\varepsilon}\chi_{C_{1}}\notag\\
&+|y^{-1}x|^{-Q-\delta}|z^{-1}x|^{-Q+\tau-\varepsilon}|wy^{-1}x|^{-Q-\beta+\varepsilon}\chi_{C_{2}}.
\end{align}
So if we take $\varepsilon \in (\beta,1)$ and integrate on $w$, we get
$$\int \chi_{(C_{1}\cup C_{2})\cap B_{1}}k_{1}(x,y,z,w)dw \lesssim|y^{-1}x|^{-Q-\beta-\delta +\varepsilon}|z^{-1}x|^{-Q+\tau-\varepsilon}.$$
Therefore, we have the desired inequality if one chooses $\varepsilon \in (\delta+\beta,\tau)$, $s=\tau-\varepsilon$ and $t=\varepsilon-\beta-\delta$.\\

The case of $C_{1}\cup C_{2}\cap B_{2}$: Notice here that 
$$|wy^{-1}z|\chi_{B_{2}}\geq c(|wy^{-1}x|-|x^{-1}z|)\chi_{B_{2}}\geq C|z^{-1}x|\chi_{B_{2}}.$$
Therefore, we have that
$$\chi_{(C_{1}\cup C_{2})\cap B_{2}}k_{1}(x,y,z,w)\lesssim |w|^{-Q-\delta}|z^{-1}x|^{-Q+\tau}|wy^{-1}x|^{-Q-\beta}.$$
Again, on $B_{2}$ we have that $|wy^{-1}x|\chi_{B_{2}}>2|z^{-1}x|\chi_{B_{2}}$. Hence,
$$\chi_{(C_{1}\cup C_{2})\cap B_{2}}k_{1}(x,y,z,w)\lesssim |w|^{-Q-\delta}|z^{-1}x|^{-Q+\tau-\varepsilon}|wy^{-1}x|^{-Q-\beta+\varepsilon}.$$
This yields
\begin{align}
\chi_{(C_{1}\cup C_{2})\cap B_{2}}k_{1}(x,y,z,w)&\lesssim |w|^{-Q-\delta}|z^{-1}x|^{-Q+\tau-\varepsilon}|y^{-1}x|^{-Q-\beta+\varepsilon}\chi_{C_{1}}\notag\\
&+|y^{-1}x|^{-Q-\delta}|z^{-1}x|^{-Q+\tau-\varepsilon}|wy^{-1}x|^{-Q-\beta+\varepsilon}\chi_{|wy^{-1}x|\lesssim |y^{-1}x|}.
\end{align}
So if we integrate over $w$ while taking $\varepsilon>\beta$, we get
$$
\int \chi_{(C_{1}\cup C_{2})\cap B_{2}}k_{1}(x,y,z,w)dw  \lesssim |y^{-1}x|^{-Q-\beta-\delta+\varepsilon}|z^{-1}x|^{-Q+\tau-\varepsilon},$$
and therefore, we have the same conclusion as in $B_{1}$.\\

Case of $(C_{1}\cup C_{2})\cap B_{3}$: Notice that on $B_{3}$ we have that $|wy^{-1}x|\chi_{B_{3}}\leq 2|z^{-1}x|\chi_{B_{3}}$ and $|wy^{-1}x|\chi_{B_{3}}>2|wy^{-1}z|\chi_{B_{3}}$, hence
\begin{align}
k_{1}\chi_{(C_{1}\cup C_{2})\cap B_{3}}&\lesssim |w|^{-Q-\delta}|wy^{-1}z|^{-Q+\tau}|wy^{-1}x|^{-Q-\beta}\notag \\
&\lesssim|w|^{-Q-\delta}|wy^{-1}z|^{-Q+\tau-\varepsilon}|wy^{-1}x|^{-Q-\beta+\varepsilon}.
\end{align}
Multiplying by $u$ and integrating with respect to $z$ yields
\begin{align}
\int \chi_{C_{1}\cup C_{2}\cap B_{3}}k_{1}(x,y,z,w)u(z)dz&\lesssim|w|^{-Q-\delta}|wy^{-1}x|^{-Q-\beta+\varepsilon}\int|wy^{-1}z|^{-Q+\tau-\varepsilon}|u|(z)dz\notag\\
&\lesssim|w|^{-Q-\delta}|wy^{-1}x|^{-Q-\beta+\varepsilon}R_{\tau-\varepsilon}|u|(w^{-1}y).
\end{align}
But on $C_{1}\cup C_{2}$ we have that
$$\max\{|y^{-1}x|,|wy^{-1}x|\}\chi_{C_{1}\cup C_{2}}\lesssim \chi_{C_{1}\cup C_{2}}|w|.$$
Thus, if we take $\varepsilon=\varepsilon_{1}+\varepsilon_{2}$, where $\varepsilon_{1}>\delta$ and $\varepsilon_{2}>\beta$, we have that
\begin{align}
\int \int \chi_{C_{1}\cup C_{2}\cap B_{3}}k_{1}(x,y,z,w)u(z)dz dw&\lesssim |y^{-1}x|^{-Q-\delta+\varepsilon_{1}}\int|wy^{-1}x|^{-Q-\beta+\varepsilon_{2}}R_{\tau-\varepsilon}|u|(w^{-1}y)dw\notag \\
&\lesssim |y^{-1}x|^{-Q-\delta+\varepsilon_{1}}R_{\tau-\beta-\varepsilon_{1}}|u|(x),
\end{align}
which finishes the proof for $s=\tau-\beta-\varepsilon_{1}$ and $t=\varepsilon_{1}-\delta$.\\

\textbf{Step 2}: Estimate on $C_{3}$:\\
In this step we will rewrite the kernel $k$ by noticing that
\begin{align}
k&=\tilde{R}_{-\delta}(w)\Big( [R_{\tau}(z^{-1}x)-R_{\tau}(z^{-1}w^{-1}y)]\tilde{R}_{-\beta}(wy^{-1}x)-[R_{\tau}(z^{-1}x)-R_{\tau}(z^{-1}y)]\tilde{R}_{-\beta}(y^{-1}x)\Big)\notag \\
&=\tilde{R}_{-\delta}(w)\Big([R_{\tau}(z^{-1}x)-R_{\tau}(y^{-1}z)][\tilde{R}_{-\beta}^{-1}(y^{-1}x)-\tilde{R}_{-\beta}^{-1}(wy^{-1}x)]\tilde{R}_{-\beta}^{-1}(wy^{-1}x)\tilde{R}_{-\beta}^{-1}(y^{-1}x)\notag\\
&+[R_{\tau}(z^{-1}y)-R_{\tau}(w^{-1}z^{-1}y)]\tilde{R}_{-\beta}(wy^{-1}x)\Big)\label{kr},
\end{align}
and we will be using the splitting induced by the $A_{i}$, $i=1,2,3$. So first we will look at $k$ on the set $C_{3}\cap(A_{1}\cup A_{2})$.\\
Notice that on $C_{3}$, $|y^{-1}x|$ and $|wy^{-1}x|$ are comparable. Moreover, 
$$\chi_{C_{3}\cap (A_{1}\cup A_{2})}|wy^{-1}z|\lesssim \chi_{C_{3}\cap (A_{1}\cup A_{2})}(|w|+|y^{-1}z|)\lesssim \chi_{C_{3}\cap (A_{1}\cup A_{2})}(|y^{-1}x|+ |y^{-1}z|)\lesssim \chi_{C_{3}\cap (A_{1}\cup A_{2})}|y^{-1}z|,$$
and
$$|wy^{-1}z|\geq c(|y^{-1}z|-|w|)\geq c(|y^{-1}z|-\frac{1}{4}|y^{-1}x|)\geq c |y^{-1}z|.$$
So again $|wy^{-1}z|$ and $|y^{-1}z|$ are comparable. Therefore, by Lemma 2.1, we have for $\varepsilon \in[0,1]$ that
\begin{align}
\chi_{C_{3}\cap (A_{1}\cup A_{2})}k&\lesssim|w|^{-Q-\delta}|y^{-1}x|^{-2Q-2\beta}\min\{|z^{-1}x|,|z^{-1}y|\}^{-Q+\tau-\varepsilon}|y^{-1}x|^{\varepsilon}|y^{-1}x|^{Q+\beta-\varepsilon}|w|^{\varepsilon}\chi_{C_{3}\cap (A_{1}\cup A_{2})}\notag\\
&+|w|^{-Q-\delta}|y^{-1}x|^{-Q-\beta}|y^{-1}z|^{-Q+\tau-\varepsilon}|w|^{\varepsilon}\chi_{C_{3}\cap (A_{1}\cup A_{2})}.
\end{align}
Taking $\varepsilon=\varepsilon_{1}+\varepsilon_{2}$, we find that
\begin{align}
\chi_{C_{3}\cap (A_{1}\cup A_{2})}k&\lesssim|w|^{-Q-\delta+\varepsilon_{1}}|z^{-1}x|^{-Q+\tau-\varepsilon}|y^{-1}x|^{-Q-\beta+\varepsilon_{2}}\notag\\
&+|w|^{-Q-\delta+\varepsilon_{1}}|z^{-1}y|^{-Q+\tau-\varepsilon}|y^{-1}x|^{-Q-\beta+\varepsilon_{2}}.
\end{align}
Once again, we choose $\varepsilon_{1}>\delta$ and $\varepsilon_{2}>\beta$, and $\varepsilon=\varepsilon_{1}+\varepsilon_{2}<\min\{1,\tau\}$, to have after integration over $w$ on $C_{3}$
$$\int \chi_{C_{3}\cap (A_{1}\cup A_{2})}k dw \lesssim\Big(|z^{-1}x|^{-Q+\tau-\varepsilon}|y^{-1}x|^{-Q+\varepsilon-\delta-\beta}+|z^{-1}y|^{-Q+\tau-\varepsilon}|y^{-1}x|^{-Q+\varepsilon-\beta-\delta}\Big).$$
Thus,
$$\int \int \int \chi_{C_{3}\cap (A_{1}\cup A_{2})}k(x,y,z,w)u(z)v(y)dwdydz\lesssim R_{s}|u|R_{t}|v|+R_{t}(|v|R_{s}|u|),$$
for $t=\varepsilon-\delta-\beta$ and $s=\tau-\varepsilon$.\\

We focus now on the last part, that is $C_{3}\cap A_{3}$.
Notice that since $|z^{-1}x|>|y^{-1}z|$ on $A_{3}$, and $|y^{-1}x|\approx |wy^{-1}x|$ on $C_{3}$, we have from $(\ref{kr})$, that
\begin{align}
\chi_{C_{3}\cap A_{3}}k\lesssim &|w|^{-Q-\delta}|y^{-1}z|^{-Q+\tau}|y^{-1}x|^{-Q-\beta}\chi_{C_{3}\cap A_{3}}\notag\\
&+|w|^{-Q-\delta}|R_{\tau}(z^{-1}y)-R_{\tau}(w^{-1}z^{-1}y)||y^{-1}x|^{-Q-\beta}\chi_{C_{3}\cap A_{3}}\notag\\
&\lesssim k_{3}+k_{4}.
\end{align}
We will estimate $k_{3}$ and $k_{4}$ separately. For $k_{3}$ we have that
\begin{align}
k_{3}&\lesssim |w|^{-Q-\delta+\varepsilon}|y^{-1}z|^{-Q+\tau}|y^{-1}x|^{-Q-\beta}|y^{-1}x|^{-\varepsilon}\chi_{C_{3}\cap A_{3}}\notag\\
&\lesssim |w|^{-Q-\delta+\varepsilon_{1}}|y^{-1}z|^{-Q+\tau-\varepsilon}|y^{-1}x|^{-Q-\beta+\varepsilon_{2}}\chi_{C_{3}\cap A_{3}},
\end{align}
which is similar to case $C_{3}\cap (A_{1}\cup A_{2})$. Hence, the desired inequality holds .\\
We move now to the term $k_{4}$. We have that
$$k_{4}\lesssim|w|^{-Q-\delta+\varepsilon_{2}}|R_{\tau}(z^{-1}y)-R_{\tau}(w^{-1}z^{-1}y)||y^{-1}x|^{-Q-\beta+\varepsilon_{2}}\chi_{C_{3}\cap A_{3}}.$$
So in the set $\{2|w|<|y^{-1}z|\}$, we have by Lemma 2.1
$$\chi_{2|w|<|y^{-1}z|}k_{4}\lesssim |w|^{-Q-\delta+\varepsilon_{1}}|z^{-1}y||y^{-1}x|^{-Q-\beta+\varepsilon_{2}}\chi_{C_{3}\cap A_{3}}\chi_{2|w|<|y^{-1}z|}$$
and we are in the same situation as $k_{3}$ above. So we treat the estimate  in $2|w|\geq|y^{-1}z|$.\\

Notice that
\begin{align}
\chi_{2|w|\geq|y^{-1}z|}k_{4}&\lesssim \chi_{C_{3}\cap A_{3}} \chi_{2|w|\geq |y^{-1}z|}|w|^{-Q-\delta+\varepsilon_{2}}|z^{-1}y|^{-Q+\tau}|y^{-1}x|^{-Q-\beta+\varepsilon_{2}}\notag\\
&+\chi_{|w^{-1}z^{-1}y|\lesssim|w|} |w|^{Q-\delta-\varepsilon_{2}+\varepsilon}|w^{-1}z^{-1}y)|^{-Q+\tau-\varepsilon}|y^{-1}x|^{-Q-\beta+\varepsilon_{2}}\chi_{C_{3}\cap A_{3}}.
\end{align}
Therefore we have first that
$$\int \chi_{2|w|\geq |y^{-1}z|}|w|^{-Q-\delta+\varepsilon_{2}}|z^{-1}y|^{-Q+\tau}|y^{-1}x|^{-Q-\beta+\varepsilon_{2}}dw \lesssim \chi_{C_{3}\cap A_{3}}|z^{-1}y|^{-Q+\tau-\delta-\varepsilon_{2}}|y^{-1}x|^{-Q-\beta+\varepsilon_{2}},$$
which provides us with the desired estimate for $\varepsilon_{1}>\delta$, $\varepsilon_{2}>\beta$ and $\varepsilon_{1}+\varepsilon_{2}=\varepsilon\in (\beta+\delta,\min\{1,\tau\})$.
On the other hand, since
$$
\int |w|^{Q-\delta-\varepsilon_{2}+\varepsilon}|w^{-1}z^{-1}y)|^{-Q+\tau-\varepsilon}|y^{-1}x|^{-Q-\beta+\varepsilon_{2}}|u(z)|dz \lesssim |w|^{Q-\delta-\varepsilon_{2}+\varepsilon}|y^{-1}x|^{-Q-\beta+\varepsilon_{2}}
R_{\tau-\varepsilon}|u|(w^{-1}y),$$
integrating again with respect to $w$ yields
$$\int |w|^{Q-\delta-\varepsilon_{2}+\varepsilon}|y^{-1}x|^{-Q-\beta+\varepsilon_{2}}
R_{\tau-\varepsilon}|u|(w^{-1}y) dw \lesssim |y^{-1}x|^{-Q-\beta+\varepsilon_{2}}R_{\tau-\tau-\varepsilon_{2}}|u|(y)$$
and in the end,
$$\int |y^{-1}x|^{-Q-\beta+\varepsilon_{2}}R_{\tau-\tau-\varepsilon_{2}}|u|(y)|v|(y)dy\lesssim R_{\varepsilon_{2}-\beta}((R_{\tau-\delta-\varepsilon_{2}}|u|)|v|)(x).$$
Which finishes the proof.

\section{Proof of Corollary 1.1}

In this case, we will write $\alpha=2k+s$ with $s\in (0,2)$, then the commutator $H_{\alpha}$ can be written as
\begin{align}
H_{\alpha}(u,v)&=(\Delta_{b})^{k}(\Delta_{b})^{\frac{s}{2}}(uv)-u(\Delta_{b})^{k}(\Delta_{b})^{\frac{s}{2}}v-v(\Delta_{b})^{k}(\Delta_{b})^{\frac{s}{2}}u\notag\\
&=(\Delta_{b})^{\frac{s}{2}}H_{k}(u,v)\notag\\
&+(\Delta_{b})^{\frac{s}{2}}(v(\Delta_{b})^{k}u)-v(\Delta_{b})^{k}(\Delta_{b})^{\frac{s}{2}}u\notag\\
&+(\Delta_{b})^{\frac{s}{2}}(u(\Delta_{b})^{k}v)-u(\Delta_{b})^{k}(\Delta_{b})^{\frac{s}{2}}v.
\end{align}
Again, in order to simplify the notations, given $m\in \mathbb{N}$, we let $\nabla_{H}^{m}$ be a combination of $\nabla_{H}$ and $\Delta_{b}$ of differentiation order $m$. Therefore, we have that
$$H_{k}(u,v)=\sum_{i=1}^{2k-1}\nabla_{H}^{i}u\nabla_{H}^{2k-i}v.$$
So we want to estimate the term $(\Delta_{b})^{\frac{s}{2}}(\nabla_{H}^{i}u\nabla_{H}^{2k-i}v)$. But first we notice that
$$(\Delta_{b})^{\frac{s}{2}}(\nabla_{H}^{i}u\nabla_{H}^{2k-i}v)=H_{s}(\nabla_{H}^{i}u,\nabla_{H}^{2k-i}v)+\nabla_{H}^{2k-i}v(\Delta_{b})^{\frac{s}{2}}\nabla_{H}^{i}u    +    \nabla_{H}^{i}u (\Delta_{b})^{\frac{s}{2}}\nabla_{H}^{2k-i}v.$$
Using Theorem 1.1, we have the existence of $s_{j,1}$, $s_{j,2}$, so that
\begin{align}
H_{s}(\nabla_{H}^{i}u,\nabla_{H}^{2k-i}v)&\lesssim \sum_{j=1}^{L}R_{s-s_{j,1}-s_{j,2}}(R_{s_{j,1}}|(-\Delta_{b})^{\frac{s}{2}}\nabla_{H}^{i}u||R_{s_{j,2}}|(-\Delta_{b})^{\frac{s}{2}}\nabla_{H}^{2k-i}v|)(x)\notag\\
&\lesssim\sum_{j=1}^{L}R_{s-s_{j,1}-s_{j,2}}(R_{s_{j,1}+\tau_{1}-s-i}|a|R_{s_{j,2}+\tau_{2}-\alpha+i}|b|)(x),
\end{align}
where we used here the fact that for $\alpha_{1}>\alpha_{2}>0$ we have
$$|R_{\alpha_{1}}(u)|=|R_{\alpha_{1}-\alpha_{2}}(\Delta_{b})^{\frac{\alpha_{2}}{2}}u)|\lesssim R_{\alpha_{1}-\alpha_{2}}|(\Delta_{b})^{\frac{\alpha_{2}}{2}}u|.$$
Therefore, we have that
\begin{align}
|(\Delta_{b})^{\frac{s}{2}}(\nabla_{H}^{i}u\nabla_{H}^{2k-i}v)|&\lesssim\Big(\sum_{j=1}^{L}R_{s-s_{j,1}-s_{j,2}}(R_{s_{j,1}+\tau_{1}-s-i}|a|R_{s_{j,2}+\tau_{2}-\alpha+i}|b|)(x)\notag\\
&+R_{\tau_{1}-s-i}|a|R_{\tau_{2}-2k+i}|b|+R_{\tau_{1}-i}|a|R_{\tau_{2}-\alpha+i}|b|\Big).
\end{align}
It follows then
\begin{align}
|(\Delta_{b})^{\frac{s}{2}}H_{k}(u,v)|&\leq C\sum_{i=1}^{2k-1}\sum_{j=1}^{L}R_{s-s_{j,1}-s_{j,2}}(R_{s_{j,1}+\tau_{1}-s-i}|a|R_{s_{j,2}+\tau_{2}-\alpha+i}|b|) \notag \\
&+R_{\tau_{1}-s-i}|a|R_{\tau_{2}-2k+i}|b|+R_{\tau_{1}-i}|a|R_{\tau_{2}-\alpha+i}|b|.
\end{align}
It remains now to estimate the terms $(\Delta_{b})^{\frac{s}{2}}(v(\Delta_{b})^{k}u)-v(\Delta_{b})^{k}(\Delta_{b})^{\frac{s}{2}}u$ and $(\Delta_{b})^{\frac{s}{2}}(u(\Delta_{b})^{k}v)-u(\Delta_{b})^{k}(\Delta_{b})^{\frac{s}{2}}v$. Notice that
$$(\Delta_{b})^{\frac{s}{2}}(v(\Delta_{b})^{k}u)-v(\Delta_{b})^{k}(\Delta_{b})^{\frac{s}{2}}u=H_{s}((\Delta_{b})^{k}u,v)+(\Delta_{b})^{k}u(\Delta_{b})^{\frac{s}{2}}v.$$
Since $\tau_{1}$ and $\tau_{2}$ are bigger than $\alpha-1=2k+s-1$, we have that $\tau_{2}>s$. We distinguish two cases.\\
\textbf{Case $\tau_{1}>2k$:}\\
In this case we have that
$$|(\Delta_{b})^{k}u|\lesssim R_{\tau_{1}-2k}|(\Delta_{b})^{\frac{\tau_{1}}{2}}u|$$
and $$|(\Delta_{b})^{\frac{s}{2}}v|\lesssim R_{\tau_{2}-s}|(\Delta_{b})^{\frac{\tau_{2}}{2}}v|$$
Combining this with Theorem 1.1, we have the desired estimate.\\
\textbf{Case $\tau_{1}\leq 2k$:}\\
We choose $\delta \in (2k-\tau_{1},1-s)$, knowing that $s<1$, we have that
\begin{align}
(\Delta_{b})^{\frac{s}{2}}(v(\Delta_{b})^{k}u)-v(\Delta_{b})^{k}(\Delta_{b})^{\frac{s}{2}}u&=(\Delta_{b})^{\frac{s}{2}}(v(\Delta_{b})^{\frac{\delta}{2}}R_{\tau_{1}-2k+\delta}(\Delta_{b})^{\frac{\tau_{1}}{2}}u)\notag\\
&\quad-v(\Delta_{b})^{\frac{s}{2}}(\Delta_{b})^{\frac{\delta}{2}}R_{\tau_{1}-2k+\delta}(\Delta_{b})^{\frac{\tau_{1}}{2}}u\notag\\
&=T_{\tau_{2},s,\delta}(V,U)
\end{align}
where $V=(\Delta_{b})^{\frac{\tau_{2}}{2}}v$ and $U=R_{\tau_{1}-2k+\delta}(\Delta_{b})^{\frac{\tau_{1}}{2}}u$. Hence, one can apply Theorem 1.2 to conclude.
\section{The case of the Heisenberg group}
In the case where $\mathbb{G}=\mathbb{H}^{n}$, and therefore $Q=2n+2$, there are two kinds of fractional powers of the sub-Laplacian. Indeed, there is the usual fractional sub-Laplacian as defined above in this manuscript, and the conformally invariant one that we will denote by $\mathcal{L}_{\alpha}$, we will refer to this last one as the geometric fractional sub-Laplacian. As in \cite{Ls}, these two operators can be distinguished in terms of their spectral Fourier multipliers (we refer the reader to \cite{than} for the definition of the Fourier transform and Fourier spectral multipliers on the Heisenberg group). Indeed, we have that the operator $(-\Delta_{b})$ corresponds to the spectral multiplier
$$ A(k,\lambda):=(2k+n)|\lambda|.$$
Therefore, $(-\Delta_{b})^{\frac{\alpha}{2}}$ corresponds to the multiplier
$$A(k,\lambda,\alpha):=((2k+n)|\lambda|)^{\frac{\alpha}{2}}.$$
On the other hand, the spectral multiplier of $\mathcal{L}_{\alpha}$ is
$$\tilde{A}(k,\lambda,\alpha):=(2|\lambda|)^{\frac{\alpha}{2}}\frac{\Gamma(\frac{2k+n}{2}+\frac{2+\alpha}{4})}{\Gamma(\frac{2k+n}{2}+\frac{2-\alpha}{4})}.$$
In particular, we have that 
$$\tilde{A}(k,\lambda,2)=A(k,\lambda,2)=A(k,\lambda),$$
and therefore $\mathcal{L}_{2}=-\Delta_{b}$.\\
As opposed to the case of a general Carnot group and the potential $R_{\alpha}$, the fundamental solution $G_{\alpha}$ of $\mathcal{L}_{\alpha}$ can be computed explicitly. Indeed, as shown in \cite{Ls} extending the result in \ref{Fol2} for the case of the sub-Laplacian, there exists a constant $c_{n,\alpha}$ such that
$$G_{\alpha}(x)=c_{n}\frac{1}{|x|^{Q-\alpha}}.$$
Also, we have a similar integral representation of $\mathcal{L}_{\alpha}$ where $\tilde{R}_{-\alpha}$ can be computed explicitely. Indeed, there exists a constant $\tilde{c}_{n,\alpha}$, such that for $u\in C^{\infty}_{0}(\mathbb{H}^{n})$ we have 
$$\mathcal{L}_{\alpha}u(x)=\tilde{c}_{n,\alpha}\int_{\mathbb{H}^{n}}\frac{u(x)-u(y)}{|y^{-1}x|^{Q+\alpha}}dy.$$
These two facts regarding the fundamental solution and the integral representation allow us then to have the same $3$-commutator result for $\mathcal{L}_{\alpha}$ that is if we define $H_{\alpha}^{\mathcal{L}}$ by
$$H_{\alpha}^{\mathcal{L}}(u,v)=\mathcal{L}_{\alpha}(uv)-u\mathcal{L}_{\alpha}v-v\mathcal{L}_{\alpha}u,$$
then we have
\begin{proposition}
Given $u,v\in C^{\infty}_{0}(\mathbb{H}^{n})$ , $0<\alpha<Q$ and $\epsilon>0$. Given $\tau_{1}$ and $\tau_{2}$ in $(\max\{0,\alpha-1\},\alpha]$ such that $\tau_{1}+\tau_{2}>\alpha$, there exists $L\in \mathbb{N}$, $s_{j,1}\in (0,\tau_{1})$, $s_{j,2}\in (0,\tau_{2})$, for $j=1,\cdots, L$, satisfying $\tau_{1}+\tau_{2}-s_{j,1}-s_{j,2}-\alpha\in [0,\epsilon)$ such that
\begin{equation}
|H_{\alpha}^{\mathcal{L}}(u,v)|(x)\lesssim \sum_{j=1}^{L}R_{\tau_{1}+\tau_{2}-s_{j,1}-s_{j,2}-\alpha}\Big(R_{s_{j,1}}|\mathcal{L}_{\tau_{1}}u|R_{s_{j,2}}|\mathcal{L}_{\tau_{2}}v|\Big)(x).
\end{equation}
\end{proposition}

{\bf Acknowledgements}
The author would like to thank Armin Schikorra for the clarifications that he provided regarding his paper \cite{Schik}.

\end{document}